\documentclass[5p,preprint]{elsarticle}

\usepackage[T1]{fontenc}

\usepackage{amsmath}
\usepackage{amssymb}
\usepackage{amsfonts}
\usepackage{amsthm}
\usepackage{tikz}
\usepackage{pgfplots}
\usepackage{graphicx,pst-all}
\usepackage{subfig}
\usepackage{enumerate}
\usepackage{caption}

\newcommand{\cC}{\mathcal{C}}

\newcommand{\cT}{\mathcal{T}}
\newcommand{\cU}{\mathcal{U}}

\newcommand{\cA}{\mathcal{A}}

\newcommand{\cL}{\mathcal{L}}

\newcommand{\cD}{\mathcal{D}}

\newcommand{\cX}{\mathcal{X}}
\newcommand{\cY}{\mathcal{Y}}

\usetikzlibrary{arrows}
\usetikzlibrary{fit}\usetikzlibrary{calc}
  \pgfdeclarelayer{background}
  \pgfsetlayers{background,main}

\theoremstyle{theorem}
\newtheorem{Prop}{Proposition}[section]

\newtheorem{Thm}[Prop]{Theorem}

\theoremstyle{definition}
\newtheorem{Def}[Prop]{Definition}
\newtheorem{Rem}[Prop]{Remark}

\newtheorem{Ex}[Prop]{Example}

\newcommand{\N}{{\mathbb{N}}}

\newcommand{\R}{{\mathbb{R}}}
\newcommand{\C}{{\mathbb{C}}}

\newcommand{\fh}{\mathfrak{h}}

\newcommand{\fL}{\mathfrak{L}}

\newcommand{\setdef}[2]{\left\{\ #1\ \left|\ \vphantom{#1} #2\ \right.\right\}}

\newcommand{\ds}[1]{{\rm \, d} #1 \,}

\DeclareMathOperator{\RE}{Re}

\DeclareMathOperator{\loc}{loc}

{\left(\begin{smallmatrix}}
{\end{smallmatrix}\right)}

{\left[\begin{smallmatrix}}
{\end{smallmatrix}\right]}

\begin{document}

\begin{frontmatter}

\title{Funnel control in the presence of infinite-dimensional internal dynamics\tnoteref{thanks}}
\tnotetext[thanks]{This work was supported by the German Research Foundation (Deutsche Forschungsgemeinschaft) via the grant BE 6263/1-1.}

\author[Paderborn]{Thomas Berger}\ead{thomas.berger@math.upb.de}
\author[Hamburg]{Marc Puche}\ead{marc.puche@uni-hamburg.de}
\author[Hamburg,Twente]{Felix L. Schwenninger}\ead{f.l.schwenninger@utwente.nl}

\address[Paderborn]{Institut f\"ur Mathematik, Universit\"at Paderborn, Warburger Str.~100, 33098~Paderborn, Germany}
\address[Hamburg]{Fachbereich Mathematik, Universit\"at Hamburg, Bundesstra{\ss}e~55, 20146~Hamburg, Germany}
\address[Twente]{Department of Applied Mathematics, University of Twente, P.O. Box 217, 7500 AE Enschede, The Netherlands}

\begin{keyword}
Adaptive control, infinite-dimensional systems, funnel control, BIBO stability.\\[5mm]
\begin{center} {\it Dedicated to the memory of Ruth F. Curtain} \end{center}
\end{keyword}

\begin{abstract}
We consider output trajectory tracking for a class of uncertain nonlinear systems whose internal dynamics may be modelled by infinite-dimensional systems which are bounded-input, bounded-output stable. We describe under which conditions these systems belong to an abstract class for which funnel control is known to be feasible. As an illustrative example, we show that for a system whose internal dynamics are modelled by a transport equation, which is not exponentially stable, we obtain prescribed performance of the tracking error.
\end{abstract}

\end{frontmatter}


%
\section{Introduction}\label{Sec:Intr}
%

We study output trajectory tracking for uncertain nonlinear systems by funnel control. As a crucial assumption, we require that the internal dynamics of the system, typically arising from a partial differential equation (PDE) in our framework,  are {\it bounded-input, bounded-output (BIBO)} stable.

Funnel control has been developed in~\cite{IlchRyan02b} for systems with relative degree one, see also the survey~\cite{IlchRyan08}. The funnel controller is a low-complexity model-free output-error feedback of high-gain type; it is an adaptive controller since the gain is adapted to the actual needed value by a time-varying (non-dynamic) adaptation scheme. Note that no asymptotic tracking is pursued, but a prescribed tracking performance is guaranteed over the whole time interval. The funnel controller proved to be the appropriate tool for tracking problems in various applications, such as temperature control of chemical reactor models~\cite{IlchTren04}, control of industrial servo-systems~\cite{Hack17} and underactuated multibody systems~\cite{BergOtto19}, speed control of wind turbine systems~\cite{Hack14,Hack15b}, DC-link power flow control~\cite{SenfPaug14}, voltage and current control of electrical circuits~\cite{BergReis14a}, oxygenation control during artificial ventilation therapy~\cite{PompAlfo14} and adaptive cruise control~\cite{BergRaue18,BergRaue19pp}.

A funnel controller for a large class of systems described by functional differential equations with arbitrary relative degree has been developed recently in~\cite{BergLe18a}. While this abstract class appears to allow for fairly general infinite-dimensional systems, cf.\ also Section~\ref{Sec:FunCon},
it is in fact not clear which types of PDE systems are encompassed. As a first result, it was shown in~\cite{BergPuch19app} that the linearized model of a moving water tank, where sloshing effects appear, belongs to the aforementioned system class. On the other hand, not even every linear, infinite-dimensional system has a well-defined (integer-valued) relative degree: In that case, results as in~\cite{IlchRyan02b,BergLe18a} cannot be applied. Instead,  the feasibility of funnel control has to be investigated directly for the (nonlinear) closed-loop system, see~\cite{ReisSeli15b} for a boundary controlled heat equation and~\cite{PuchReis19pp} for a general class of boundary control systems.

The present paper is devoted to systems which have a relative degree, but in the presence of internal dynamics that are modelled by a PDE system. Motivated by the observation that several relevant systems of the aforementioned form belong to the class introduced in~\cite{BergLe18a}, we develop a general system class containing PDE models for which funnel control is feasible; this result is presented in Section~\ref{Sec:operators}. We show that the class of systems for which a Byrnes-Isidori form exists, see~\cite{IlchSeli16}, is contained in this new system class. As an example, we consider a system internally driven by a transport equation and illustrate the funnel controller by a simulation in  Section~\ref{Sec:Sim}. Some conclusions are given in Section~\ref{Sec:Concl}.

\subsection{Nomenclature and basic concepts}\label{Ssec:Nomencl}

Throughout this article, we use the following notation: $\N$ denotes the natural numbers, $\N_0 = \N \cup\{0\}$, and $\R_{\ge 0} =[0,\infty)$. We use the notation $\C_{\omega}=\setdef{\lambda\in\C}{\RE \lambda>\omega}$ for $\omega\in\R$. With $L^p(I;\R^n)$ we denote the Lebesgue space of all measurable and $p$th power integrable functions $f:I\to\R^n$, where $I\subseteq\R$ is an interval and $p\in[1,\infty)$; $L^\infty(I;\R^n)$ denotes the Lebesgue space of all measurable and essentially bounded functions $f:I\to\R^n$. We write $\|\cdot\|_\infty$ for $\|\cdot\|_{L^\infty(\R_{\ge 0};\R^n)}$. By $L^\infty_{\loc}(I;\R^n)$ we denote the set of measurable and locally essentially bounded functions $f:I\to\R^n$ and by $W^{k,p}(I;\R^n)$, $k\in\N_0$, the Sobolev space of $k$-times weakly differentiable functions $f:I\to\R^n$ such that $f,\ldots, f^{(k)}\in L^p(I;\R^n)$. For an open set $V\subseteq\R^m$ we denote by $\cC^k(V;\R^n)$ the set of $k$-times continuously differentiable functions $f:V\to\R^n$, $k\in\N_0\cup\{\infty\}$ where $\cC(V;\R^n):=\cC^0(V;\R^n)$. 
The set of all real-valued Borel measures with bounded total variation is denoted by ${\rm M}(\R_{\ge0})$ and the total variation by $\|f\|_{{\rm M}(\R_{\ge0})}$ for $f\in{\rm M}(\R_{\ge0})$; we refer to the textbook~\cite{Graf14} for more details. By $\cL(\cX;\cY)$, where $\cX, \cY$ are Hilbert spaces, we denote the set of all bounded linear operators $\cA:\cX\to\cY$.

Let $\cX$ be a real Hilbert space and recall that a {\it $\cC_0$-semigroup}  $(T(t))_{t\ge0}$ on $\cX$ is a $\mathcal{L}(\cX;\cX)$-valued map satisfying $T(0)=I_{\cX}$ and $T(t+s)=T(t)T(s)$, $s,t\geq0$, where $I_{\cX}$ denotes the identity operator, and $t\mapsto T(t)x$ is continuous for every $x\in \cX$. $\cC_0$-semigroups are characterized by their generator~$A$, which is a, not necessarily bounded, operator on~$\cX$.

Furthermore, recall the space $\cX_{-1}$, see e.g.~\cite[Sec.~2.10]{TucsWeis09}, which should be thought of as an abstract Sobolev space with negative index\footnote{This space is sometimes referred to as {\it rigged Hilbert space}.}. If $A:\cD(A)\subseteq \cX\to\cX$ is a densely defined operator with $\rho(A)\neq \emptyset$, where $\rho(A)$ denotes the resolvent set of~$A$, then for any $\beta\in\rho(A)$ we denote by $\cX_{-1}$ the completion of $\cX$ with respect to the norm
\[
    \|x\|_{\cX_{-1}}=\|(\beta I-A)^{-1}x\|_\cX,\quad x\in \cX.
\]
Then the norms generated as above for different $\beta\in\rho(A)$ are equivalent and, in particular, $\cX_{-1}$ is independent of the choice of~$\beta$. If $A$ generates a $\cC_0$-semigroup $(T(t))_{t\ge 0}$ in $\cX$, then the latter has a unique extension to a semigroup $(T_{-1}(t))_{t\ge 0}$ in~$\cX_{-1}$, which is given by
\[
    T_{-1}(t) = (\beta I - A_{-1}) T(t),\quad t\ge 0,
\]
where $(\beta I - A_{-1}) \in \cL(\cX;\cX_{-1})$ is a surjective isometry. Therefore, $A_{-1}$ is the generator of the semigroup $(T_{-1}(t))_{t\ge 0}$.

The notion of {\it admissible} operators is well-known in infinite-dimensional linear systems theory with unbounded control and observation operators, as present in boundary control, see e.g.~\cite{TucsWeis09}, and is motivated by interpreting a PDE on a larger space in order to define solutions. Let $\cU, \cX, \cY$ be real Hilbert spaces and $A$ as above such that it generates a $\cC_0$-semigroup $(T(t))_{t\ge 0}$ on $\cX$. Then we recall that $B\in \mathcal{L}(\cU;\cX_{-1})$ is a {\it $L^p$-admissible}  control operator (for $(T(t))_{t\ge 0}$), with $p\in [1,\infty]$, if for all $t\ge 0$ and all $u\in L^{p}([0,t];\cU)$ we have
$$\Phi_{t}u :=\int_{0}^{t} {T}_{-1}(t-s)Bu(s)\ds{s} \in \cX.$$
By a closed graph theorem argument this property implies that, for any $t\ge 0$, the operator $\Phi_{t}$ is bounded from  $L^{p}([0,t];\cU)$ to $\cX$.

An operator $C\in\mathcal{L}(\cD(A);\cY)$ is called {\it $L^{p}$-admissible} observation operator (for $(T(t))_{t\ge 0}$), if for some (and hence all) $t\ge 0$  the mapping
$$\Psi_{t}:\cD(A)\to L^{p}([0,t],\cY),\ x\mapsto C T(\cdot)x$$
can be extended to a bounded operator from $\cX$ to $L^{p}([0,t],\cY)$ --- this extension will again be denoted by $\Psi_{t}$.

Both admissibility notions are combined in the stronger concept of {\it well-posedness}: Let $({A},B,C)$ represent a system where  ${A}$ is the generator of a $\cC_0$-semigroup, $B$ is a $L^{2}$-admissible control operator and $C$ is a $L^{2}$-admissible observation operator in the sense described above. If for some $\omega\in\R$ the transfer function $H:\C_\omega\to\mathcal{L}(\mathcal{U},\mathcal{Y})$, which is uniquely determined (up to a constant) by
  \[
    \frac{1}{s_2-s_1}(H(s_1) - H(s_2)) = C\big( (s_1I-A)^{-1}(s_2I-A)^{-1}\big) B
  \]
  for all $s_1,s_2\in\C_\omega,\ s_1\neq s_2$, exists and is proper, that is $\sup_{s\in\C_\omega} \|H(s)\| <\infty$, then we say that $(A,B,C)$ is {\it well-posed}. We remark that well-posedness is usually defined differently, but equivalently, see~\cite{CurtWeis89}. If $\lim_{\RE s\to\infty} H(s) v$ exists for any $v\in \mathcal{U}$, then the system $(A,B,C)$ is called {\it regular}.

\subsection{System class}\label{Ssec:SysClass}

In the remainder of the present paper we consider abstract differential equations of the form
\begin{equation}\label{eq:nonlSys}
\begin{aligned}
y^{(r)}(t)&= f\big(d(t), T(y,\dot{y},\ldots,y^{(r-1)})(t)\big)\\
&\quad + \Gamma\big(d(t), T(y,\dot{y},\ldots,y^{(r-1)})(t)\big)\, u(t)\\
y|_{[-h,0]} &= y^0\in W^{r-1,\infty}([-h,0];\R^m),
\end{aligned}
\end{equation}
where $h\ge 0$ is the ``memory'' of the system\footnote{Here, ``$h=0$'' means that the initial values $y(0), \dot y(0), \ldots$, $y^{(r-1)}(0)$ are prescribed.}, $r\in\N$ is the relative degree, and
\begin{itemize}
\item[(N1)] the disturbance satisfies $d\in L^{\infty}(\R_{\ge 0};\R^p)$, $p\in \N$;
\item[(N2)] $f\in \mathcal{C}(\R^p\times \R^q; \R^m),\ q\in \N$;
\item[(N3)] the high-frequency gain matrix function $\Gamma\in\cC(\R^p\times\R^q;\R^{m\times m})$ satisfies $\Gamma(d,\eta) + \Gamma(d,\eta)^\top > 0$ for all $(d,\eta)\in\R^p\times\R^q$;
\item[(N4)] $T:\mathcal{C}([-h,\infty);\R^{rm})\rightarrow L_{\rm loc}^{\infty}(\R_{\ge 0};\R^q)$ is an operator with the following properties:
\begin{itemize}
\item[a)] $T$ maps bounded trajectories to bounded trajectories, i.e, for all $c_1>0$, there exists $c_2>0$ such that for all $\zeta\in \mathcal{C}([-h,\infty);\R^{rm})$,
\[\sup\limits_{t\in [-h,\infty)}\|\zeta(t)\|\le c_1\Rightarrow \sup\limits_{t\ge 0}\|T(\zeta)(t)\|\le c_2,\]
\item[b)] $T$ is causal, i.e, for all $t\ge 0$ and all $\zeta,\xi\in\mathcal{C}([-h,\infty);\R^{rm})$,
\[\zeta|_{[-h,t)}=\xi|_{[-h,t)}\Rightarrow T(\zeta)|_{[0,t)}\overset{\rm a.e.}{=}T(\xi)|_{[0,t)}.\]
\item[c)] $T$ is locally Lipschitz continuous in the following sense: for all $t\ge 0$ and all $\xi\in\cC([-h,t];\R^{rm})$ there exist $\tau, \delta, c>0$ such that, for all $\zeta_1,\zeta_2\in \mathcal{C}([-h,\infty);\R^{rm})$ with $\zeta_i|_{[-h,t]}=\xi$ and $\|\zeta_i(s)-\xi(t)\|<\delta$ for all $s\in[t,t+\tau]$ and $i=1,2$, we have
    \begin{multline*}
        \left\|\left(T(\zeta_1)-T(\zeta_2)\right)|_{[t,t+\tau]} \right\|_{\infty} \\
        \le c\left\|(\zeta_1-\zeta_2)|_{[t,t+\tau]}\right\|_{\infty}.
    \end{multline*}
\end{itemize}
\end{itemize}

In~\cite{BergLe18a, HackHopf13, IlchRyan09, IlchRyan02b, IlchRyan07} it is shown that the class of systems~\eqref{eq:nonlSys} encompasses linear and nonlinear systems with strict relative degree~$r$ and BIBO stable internal dynamics. The operator~$T$ allows for infinite-dimensional (linear) systems, systems with hysteretic effects or nonlinear delay elements, and combinations thereof. Note that~$T$ is typically the solution operator corresponding to a (partial) differential equation which describes the internal dynamics of the system. The linear infinite-dimensional systems that are considered in~\cite{IlchRyan02b, IlchRyan07} are in a special Byrnes-Isidori form that is discussed in detail in~\cite{IlchSeli16}. While the internal dynamics in these systems is allowed to correspond to a strongly continuous semigroup, all other operators are assumed to be bounded and to satisfy additional restrictive conditions. In contrast to this, in the present paper we consider nonlinear equations which, in particular, involve unbounded operators. This complements and generalizes the findings in~\cite{BergPuch19app}.

\subsection{Control objective}\label{Ssec:ContrObj}

The objective is to design a derivative output error feedback of the form
\[u(t) = G\big(t, e(t), \dot e(t), \ldots, e^{(r-1)}(t)\big),\]
where $y_{\rm ref}\in W^{r,\infty}(\R_{\ge 0};\R^m)$ is a reference signal, which applied to~\eqref{eq:nonlSys} results in a closed-loop system where the tracking error $e(t)=y(t)-y_{\rm ref}(t)$ evolves within a prescribed performance funnel
\begin{equation}
\mathcal{F}_{\varphi} := \setdef{(t,e)\in\R_{\ge 0} \times\R^m}{\varphi(t) \|e\| < 1},\label{eq:perf_funnel}
\end{equation}
which is determined by a function~$\varphi$ belonging to
\[
 \Phi_r :=
\setdef{\!
\varphi\in  \cC^r(\R_{\ge 0};\R)\!
}{\!\!\!\!
\begin{array}{l}
\text{ $\varphi, \dot \varphi, \ldots, \varphi^{(r)}$ are bounded,}\\
\text{ $\varphi (\tau)>0$ for all $\tau>0$,}\\
 \text{ and }  \liminf_{\tau\rightarrow \infty} \varphi(\tau) > 0
\end{array}
\!\!\!\!\!}.
\]
Furthermore, all signals $u, e, \dot e, \ldots, e^{(r-1)}$ should remain bounded.

\begin{figure}[h]
\hspace{6mm}\includegraphics[width=8cm, trim=170 550 230 120, clip]{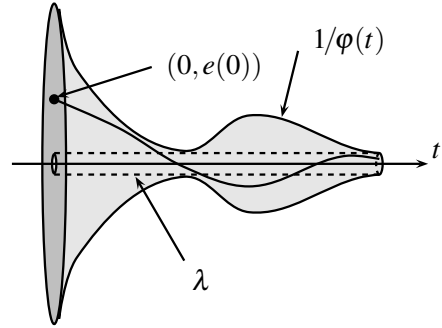}
 \vspace*{-2mm}
 \caption{Error evolution in a funnel $\mathcal F_{\varphi}$ with boundary~$1\!/\!\varphi(t)$.}
 \label{Fig:funnel}
 \end{figure}

The funnel boundary is given by~$1\!/\!\varphi$, see Fig.~\ref{Fig:funnel}. The case $\varphi(0)=0$ is explicitly allowed and puts no restriction on the initial value since $\varphi(0) \|e(0)\| < 1$; in this case the funnel boundary $1/\varphi$ has a pole at $t=0$.

An important property is that each performance funnel $\mathcal{F}_{\varphi}$ with $\varphi\in\Phi_r$ is bounded away from zero, because boundedness of~$\varphi$ implies existence of $\lambda>0$ such that $1\!/\!\varphi(t)\geq\lambda$ for all $t > 0$. The funnel boundary is not necessarily monotonically decreasing
and there are situations, like in the presence of periodic disturbances, where widening the funnel over some later time interval might be beneficial.
For typical choices of funnel boundaries see e.g.~\cite[Sec.~3.2]{Ilch13}.

\section{Funnel control}\label{Sec:FunCon}
%

It was shown in~\cite{BergLe18a} that the funnel controller
\begin{equation}\label{eq:fun-con}
\boxed{\begin{aligned}
u(t) &= - k_{r-1}(t) \, e_{r-1}(t),\\
e_0(t)&=e(t) = y(t) - y_{\rm ref}(t),\\
e_1(t)&=\dot{e}_0(t)+k_0(t)\,e_0(t),\\
e_2(t)&=\dot{e}_1(t)+k_1(t)\,e_1(t),\\
& \ \ \vdots \\
e_{r-1}(t)&=\dot{e}_{r-2}(t)+k_{r-2}(t)\,e_{r-2}(t),\\
k_i(t)&=\frac{1}{1-\varphi_i(t)^2\|e_i(t)\|^2},\quad i=0,\dots,r-1,
\end{aligned}
}
\end{equation}
where
\begin{equation}\label{eq:cond-phi}
    \varphi_0\in\Phi_r, \varphi_1\in\Phi_{r-1},\ldots,\varphi_{r-1} \in\Phi_1,
\end{equation}
achieves the control objective described in Section~\ref{Ssec:ContrObj} for any system which belongs to the class~\eqref{eq:nonlSys}. We stress that while the derivatives $\dot e_0,\ldots, \dot e_{r-2}$ appear in~\eqref{eq:fun-con}, they only serve as short-hand notations and may be resolved in terms of the tracking error, the funnel functions and the derivatives of these, cf.~\cite[Rem.~2.1]{BergLe18a}.

The existence of solutions of the initial value problem resulting from the application of the funnel controller~\eqref{eq:fun-con} to a system~\eqref{eq:nonlSys} must be treated carefully. By a \emph{solution} of~\eqref{eq:fun-con},~\eqref{eq:nonlSys} on $[-h,\omega)$ we mean a function $y\in\cC^{r-1}([-h,\omega);\R^m)$, $\omega\in(0,\infty]$, with $y|_{[-h,0]} = y^0$ such that $y^{(r-1)}|_{[0,\omega)}$ is weakly differentiable and satisfies the differential equation in~\eqref{eq:nonlSys} with~$u$ defined in~\eqref{eq:fun-con} for almost all $t\in[0,\omega)$; $y$ is called \emph{maximal}, if it has no right extension that is also a solution. Existence of solutions of functional differential equations has been investigated in~\cite{IlchRyan02b} for instance.

The following result is from~\cite{BergLe18a}. Note that in~\cite{BergLe18a} a slightly stronger version of conditions~(N3) and~(N4)~c) is used. However, the proof does not change; in particular, regarding~(N4)~c), the existence part of the proof in~\cite{BergLe18a} relies on a result from~\cite{IlchRyan09} where the version from the present paper is used.

\begin{Thm}\label{Thm:FunCon-nonl}
Consider a~system~\eqref{eq:nonlSys} with properties (N1)--(N4) for some $r\in\N$ and $h\ge 0$. Let $y_{\rm ref}\in W^{r,\infty}(\R_{\ge 0};\R^m)$, $\varphi_0,\ldots,\varphi_{r-1}$ as in~\eqref{eq:cond-phi} and $y^0\in W^{r-1,\infty}([-h,0];\R^m)$ be an initial condition such that $e_0,\ldots,e_{r-1}$ defined in~\eqref{eq:fun-con} satisfy
\[
    \varphi_i(0) \|e_i(0)\| < 1\quad \text{for}\ i=0,\ldots,r-1.
\]
Then the funnel controller~\eqref{eq:fun-con} applied to~\eqref{eq:nonlSys} yields an initial-value problem which has a solution, and every solution can be extended to a maximal solution $y:\left[-h,\omega\right)\rightarrow \R^m$, $\omega\in(0,\infty]$, which has the following properties:
\begin{enumerate}
\item The solution is global, i.e., $\omega=\infty$.
\item The input $u:\R_{\ge0}\to\R^m$, the gain functions $k_0, \ldots, k_{r-1}:\R_{\ge0}\to\R$ and $y, \dot y, \ldots, y^{(r-1)}:\R_{\ge0}\to\R^m$ are bounded.
\item The functions $e_0,\ldots,e_{r-1}:\R_{\ge 0}\to\R^m$ evolve in their respective performance funnels and are uniformly bounded away from the funnel boundaries in the sense
\begin{multline*}
  \forall\, i=0,\ldots,r-1\ \exists\, \varepsilon_i>0\ \forall\, t>0:\\
  \|e_i(t)\|\le \varphi_i(t)^{-1} - \varepsilon_i.
\end{multline*}
\end{enumerate}
\end{Thm}

While the class of functional differential equations~\eqref{eq:nonlSys} appears to be rather general and funnel control is feasible for these systems by Theorem~\ref{Thm:FunCon-nonl}, it is not clear exactly which kind of systems that contain PDEs are encompassed by the class~\eqref{eq:nonlSys}. The operator~$T$, which describes the internal dynamics, is able to model a broad class of PDE systems, as we will show in the following example which motivates the introduction of the operator class in Section~\ref{Sec:operators}.

\begin{Ex}\label{Ex:TranspEqn}
Consider the following system whose internal dynamics are described by a transport equation, that is
\begin{equation}\label{eq:sys-sim}
\begin{aligned}
\dot{y}(t)={}&z(t,0)+\gamma u(t)\\
\frac{\partial z}{\partial t}(t,\xi)={}&c\frac{\partial z}{\partial \xi}(t,\xi)+\fh(\xi)y(t),\\
z(0,\xi)={}&0,
\end{aligned}
\end{equation}
for $(t,\xi)\in(0,\infty)\times[0,\infty)$, where $c>0$ and $\fh\in {\rm M}(\R_{\ge0})$ is a Borel measure of bounded total variation. It is well-known that the second and third equations in~\eqref{eq:sys-sim} constitute a regular well-posed linear system $(A,B,C)$ on $X=L^{2}(\R_{\ge0};\R)$, the so-called {\it shift-realization} of the Laplace transform~$\fL(\fh)$, see e.g.~\cite{Helt76,Yama81}. More precisely, the PDE is then considered on the abstract Sobolev space $X_{-1}$ to appropriately interpret the term~$\fh(\xi)y(t)$ and the solutions are {\it mild solutions}\footnote{See e.g.~\cite{TucsWeis09} for a definition of the mild solution.} in general.

Also note that the generated (left-) shift-semigroup is not exponentially stable. In particular, the Laplace transform~$\fL(\fh)$ of the measure~$\fh$ is defined on the closed right half-plane and bounded analytic on this domain. Moreover, the impulse response of the PDE equals~$\fh$. More precisely, for sufficiently smooth~$y$ we have the representation
\[
    z(t,0)=(\fh\ast y)(t)=\int_{0}^{t}y(t-s)\,{\rm d}\fh(s).
\]
Therefore, the first equation in \eqref{eq:sys-sim} formally reads
\begin{equation}\label{eq:TE-new}
    \dot{y}(t)=(\fh\ast y)(t)+\gamma u,
\end{equation}
which is an integral-differential Volterra equation. Also note that for the following simple cases
\begin{itemize}
\item $\fh=\delta_{0}$, we obtain a finite-dimensional linear system: $$\dot{y}(t)=y(t)+\gamma u(t);$$
\item $\fh=\delta_{t_{0}}$, $t_0 >0$, we obtain a delay differential equation:
$$\dot{y}(t)=\begin{cases} y(t-t_{0})+\gamma u(t), & t\ge t_0,\\ \gamma u(t), & 0\le t < t_0.\end{cases}$$
\end{itemize}
Another typical case is that $\fh(\xi)=f(\xi){\rm d}\xi$ with $f\in L^{1}(\R_{\ge0};\R)$, i.e.,~$\fh$ is represented by its $L^{1}$-density with respect to the Lebesgue measure. If additionally $f\in L^{2}(\R_{\ge0};\R)$, then the input operator~$B=\fh$ of the PDE is bounded.

We may now observe that~\eqref{eq:TE-new} belongs to the system class~\eqref{eq:nonlSys}, if we define the operator
\[
    T(y) := \fh\ast y,\quad y\in\cC(\R_{\ge 0};\R).
\]
As~$\fh$ has bounded total variation, it follows that~$T$ is a bounded operator from $\cC(\R_{\ge0};\R)\cap L^{\infty}(\R_{\ge0};\R)$ to $L^{\infty}(\R_{\ge0};\R)$ and hence it is straightforward to check that~$T$ satisfies condition~(N4).
\end{Ex}

%
\section{A class of operators for funnel control}\label{Sec:operators}
%

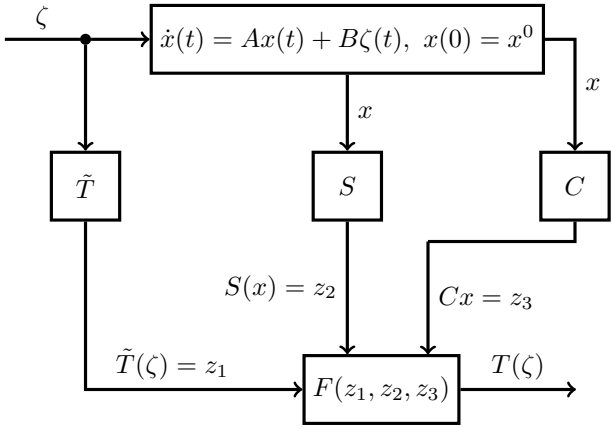
\begin{figure}[h!tb]
\begin{center}
\resizebox{9cm}{!}{
\begin{tikzpicture}[very thick,node distance = 10ex, box/.style={fill=white,rectangle, draw=black}, blackdot/.style={inner sep = 0, minimum size=3pt,shape=circle,fill,draw=black},plus/.style={fill=white,circle,inner
sep = 0,thick,draw},metabox/.style={inner sep = 0ex,rectangle,draw,dotted,fill=gray!20!white}]

    \node (PDE)     [box,minimum size=6ex]  {$\dot x(t) = Ax(t) + B\zeta(t),\ x(0) = x^0$};
    \node (fork1)   [blackdot,left of = PDE, xshift = -13ex, yshift = 0ex] {};
    \node (in)      [left of = fork1, xshift = 2ex, yshift = 0ex] {};
    \node (Ttilde)  [box,minimum size=6ex,below of = fork1, yshift = -3ex]  {$\tilde T$};
    \node (S)  [box,minimum size=6ex,below of = PDE, yshift = -3ex] {$S$};
    \node (C)  [box,minimum size=6ex,right of = S, xshift = 10ex] {$C$};
    \node (Fin1)  [below of = S, yshift = -5.8ex]  {};
    \node (Fin2)  [right of = Fin1, xshift = -2ex, yshift=0.9ex] {};
    \node (Cin)  [above of = Fin2, yshift = 0ex]  {};
    \node (F)  [box,minimum size=6ex,below of = S, yshift = -8ex,xshift=3ex]  {$F(z_1,z_2,z_3)$};
    \node (out)  [right of = F, xshift = 8ex]  {};
    \node (out2)  [right of = out, xshift = -5ex]  {};

    \draw[-] (in) -- (fork1) node[midway,above] {$\zeta$};
    \draw[->] (fork1) -- (PDE) node[midway,above] {};
    \draw[->] (fork1) -- (Ttilde) node[midway,above] {};
    \draw[->] (Ttilde) |- (F) node[pos=0.7,above] {$\tilde T(\zeta) = z_1$};
    \draw[->] (PDE) -- (S) node[midway,right] {$x$};
    \draw[->] (PDE) -| (C) node[pos=0.7,right] {$x$};
    \draw[->] (S) -- (Fin1) node[midway,left] {$S(x)=z_2$};
    \draw[-] (C) |- (Cin.west) node[midway,right] {};
    \draw[->] (Cin.west) -- (Fin2.west) node[midway,right] {$Cx=z_3$};
    \draw[->] (F) -- (out) node[midway,above] {$T(\zeta)$};
  \end{tikzpicture}
}
\end{center}
  \caption{Structure of an operator $T\in\cT_h^{\ell,q}$.} \label{Fig:OpStruc}
\end{figure}

Motivated by Example~\ref{Ex:TranspEqn}, in this section we develop a description for a class of operators~$T$ which include certain linear PDEs and satisfy condition~(N4). The aforementioned PDEs may either be coupled with a nonlinear observation operator which satisfies a certain growth bound, or it may be coupled with a linear observation operator which is possibly unbounded, but with respect to which the system is regular well-posed. In both cases we additionally require that the overall system is BIBO stable. For the linear observation operator, this is true if, for instance, the inverse Laplace transform of the corresponding transfer function defines a Borel measure with bounded total variation. This structure is illustrated in Fig.~\ref{Fig:OpStruc}.

We give a precise definition of the operator class in the following.

\begin{Def}\label{Def:opT}
Let $h\ge 0$ and $\ell,q\in\N$. Then $\cT_h^{\ell,q}$ is defined as the set of all operators
\[T:\cC([-h,\infty);\R^\ell)\to L^\infty_{\loc}(\R_{\ge 0};\R^q)\]
 which, for any $\zeta\in\cC([-h,\infty);\R^\ell)$, are given by
\[
    T(\zeta)(t) = F\big( \tilde T(\zeta)(t), S(x)(t), (Cx)(t)\big),\quad t\ge 0,
\]
where $x$, for some $x^0 \in \cD(A)$, is the mild solution of the PDE
\begin{equation}\label{eq:PDE}
   \dot x(t) = Ax(t) + B\zeta(t),\quad x(0) = x^0,
\end{equation}
where
\begin{itemize}
  \item[(P1)] $A$ generates a bounded $\cC_0$-semigroup in a real Hilbert space $X$ and $B\in\cL(\R^\ell;X_{-1})$,  $C\in\cL(\cD(A);\R^{q_3})$ are operators such that $(A,B,C)$ is a regular well-posed linear system which additionally is BIBO stable, i.e., the operator
$$L^{\infty}((0,\infty);\R^{\ell})\to L^{\infty}((0,\infty);\R^{q_{3}}),\ f\mapsto \fL^{-1}(H)\ast f $$
is bounded, where $H:\C_0\to\C^{q_3\times\ell}$
 denotes the transfer function of $(A,B,C)$.
  \item[(P2)] $F\in\cC^1(\R^{q_1}\times\R^{q_2}\times\R^{q_3};\R^q)$;
  \item[(P3)] $\tilde T:\cC([-h,\infty);\R^\ell)\to L^\infty_{\loc}(\R_{\ge 0};\R^{q_1})$ satisfies condition~(N4) in Section~\ref{Ssec:SysClass} with $\ell = rm$;
  \item[(P4)] $S:X\to\R^{q_2}$ is a (possibly nonlinear) operator which satisfies that for all $x\in X$ and all $\rho>0$ there exists $L>0$ such that for all $x_1, x_2\in X$ with $\|x_i - x\|_X < \rho$, $i=1,2$, we have
      \[
        \|S(x_1) - S(x_2)\| \le L \|S(x_1 - x_2)\|.
      \]
   Furthermore,~$S$ is such that~\eqref{eq:PDE} is BIBO stable w.r.t.~$S$, i.e., there exists $\gamma\in \cC^{1}(\R_{\ge0};\R)$ such that for all $\zeta\in\cC([-h,\infty);\R^\ell)$ the mild solution of~\eqref{eq:PDE} satisfies
      \[
       \forall\, t\ge0:\quad \|S\big(x(t)\big)\| \le \gamma (\|\zeta|_{[-h,t]}\|_\infty);
      \]
\end{itemize}
\end{Def}

\begin{Rem}\
\begin{enumerate}
\item We note that any operator~$T$ as given in Definition~\ref{Def:opT} with the properties~(P1)--(P1) is indeed well-defined from $\cC([-h,\infty);\R^\ell)$ to $L^\infty_{\loc}(\R_{\ge 0};\R^q)$.

  \item We emphasize that the assumption of BIBO stability of~\eqref{eq:PDE} as in~(P4) is quite weak. Provided that~$S$ is sufficiently nice, then a sufficient condition for this is {\it input-to-state stability}~\cite{Sont89a}. This concept was studied extensively for nonlinear systems, see~\cite{Sont08}, and for systems containing  PDEs it is investigated in~\cite{JacoNabi18,MiroWirt18}. However, the state of an input-to-state stable system converges to zero whenever the input is zero, which is not required for BIBO stable systems considered here.

  \item Note that the assumption of BIBO stability in~(P1) essentially reduces to showing that the inverse Laplace transform $\fh_{ij} = \fL^{-1}(H_{ij})$ is a Borel measure on $\R_{\ge 0}$ with bounded total variation for all $i=1,\ldots,q_3$ and $j=1,\ldots,\ell$, i.e., $\fh_{ij} \in {\rm M}(\R_{\ge 0})$. Recall that there exist bounded, shift-invariant operators on $L^{\infty}((0,\infty);\R)$ defined as the convolution with a tempered distribution, which is not contained in ${\rm M}(\R_{\ge 0})$, see \cite[Sec.\ 2.5.4]{Graf14}.

\end{enumerate}
\end{Rem}

In the following main result we show that any operator which belongs to the class $\cT_h^{\ell,q}$ satisfies the condition~(N4) in Section~\ref{Ssec:SysClass}.

\begin{Thm}\label{mainthm}
Any $T\in\cT_h^{\ell,q}$ satisfies condition~(N4) in Section~\ref{Ssec:SysClass}.
\end{Thm}
\begin{proof}
\emph{Step 1}: We show property~(N4)~a). To this end, observe that by continuity of~$F$ it suffices to show this for the maps $\zeta\mapsto \tilde T(\zeta)$, $\zeta\mapsto S(x)$ and $\zeta\mapsto Cx$; recall that~$x$ as in~\eqref{eq:PDE} depends on~$\zeta$. By~(P3), $\tilde T$ satisfies~(N4)~a) and by~(P4) we have
\begin{align*}
    \|S(x(t))\| \le \gamma (\|\zeta\|_\infty)
\end{align*}
for all $t\ge 0$ and all bounded $\zeta\in\cC([-h,\infty);\R^\ell)$. It remains to show that~$Cx$ is bounded. By~(P1) the system $(A,B,C)$ is regular and well-posed, from which it follows by the variation of constants formula, see e.g.~\cite{TucsWeis14}, that
\[Cx(\cdot)=CT_A(\cdot)x_0+(\fh\ast\zeta)(\cdot),\]
where $(T_A(t))_{t\geq0}$ is the $\cC_0$-semigroup generated by~$A$ and $\fh = \fL^{-1}(H)$ is the inverse Laplace transform of the transfer function~$H:\C_0\to\C^{q_3\times\ell}$. By Assumption~(P1) there exists $C_\fh>0$ such that $\|\fh\ast \zeta\|_{\infty}\leq C_{\fh}\|\zeta\|_{\infty}$ and
  thus, for all $t\ge 0$,
\begin{align*}
\|Cx(t)\|\leq&\ \|CT_A(t)x_0\|+\|(\fh\ast\zeta)(t)\|\\
\leq&\ \|C\|_{\cL(\cD(A);\R^{q_3})}\|AT_A(t)x_0\|+C_{\fh}\|\zeta\|_{\infty}\\
=&\ \|C\|_{\cL(\cD(A);\R^{q_3})}\|T_A(t)Ax_0\|+C_{\fh}\|\zeta\|_{\infty}\\
\leq&\ \|C\|_{\cL(\cD(A);\R^{q_3})}\|T_A(t)\|_{\cL(X)}\|Ax_0\|_X+C_{\fh}\|\zeta\|_{\infty}\\
\leq&\ M\|C\|_{\cL(\cD(A);\R^{q_3})}\|Ax_0\|_X+C_{\fh}\|\zeta\|_{\infty},
\end{align*}
where we have used that $x_0\in\cD(A)$ and $(T_A(t))_{t\geq0}$ is bounded, that is, $M=\sup_{t\geq0}\|T_A(t)\|_{\cL(X;X)}<\infty$. Thus,
\[\|Cx(\cdot)\|_\infty\leq M \|C\|_{\cL(\cD(A);\R^{q_3})}\|Ax_0\|_X+C_{\fh}\|\zeta\|_{\infty}.\]

\emph{Step 2}: We show property~(N4)~b). This is a straightforward consequence of the definition of~$\tilde T$.

\emph{Step 3}: We show property~(N4)~c). Fix $t\ge 0$ and $\xi\in\cC([-h,t];\R^\ell)$. Let $\tilde \tau, \tilde\delta, \tilde c$ be the constants given by property~(N4)~c) of~$\tilde T$. Set $\tau:=\tilde \tau$ and $\delta := \tilde \delta$. Further let $\zeta_i\in\cC([-h,\infty);\R^\ell)$ with $\zeta_i|_{[-h,t]}=\xi$ and $\|\zeta_i(s)-\xi(t)\|<\delta$ for all $s\in[t,t+\tau]$ and $i=1,2$. Let~$x_i$ denote the mild solution of~\eqref{eq:PDE} corresponding to~$\zeta_i$ for $i=1,2$. Then, by linearity, $x_1-x_2$ is the mild solution corresponding to~$\zeta_1-\zeta_2$. 
Let~$\tilde x$ denote the mild solution of~\eqref{eq:PDE} corresponding to~$\tilde \xi$ defined by $\tilde \xi|_{[-h,t]} = \xi$ and $\tilde \xi|_{[t,\infty)} \equiv \xi(t)$. Then, since by well-posedness of $(A,B,C)$ the operator  $B$ is $L^2$-admissible, we have for all $s\in [t,t+\tau]$ that
\begin{align*}
    \|x_i(s) - \tilde x(t)\|_X &\le \|\Phi_{t+\tau}\big((\zeta_i-\xi(t))|_{[t,s]}\big)\|_X \\
    &< \delta \|\Phi_{t+\tau}\|.
\end{align*}
Now let~$L$ be the constant given by~(P4) for $x= \tilde x(t)$ and $\rho = \delta \|\Phi_{t+\tau}\|$, and further set
\[
   L_2 := L \cdot \sup_{s\in[0,2\delta]} |\gamma'(s)|.
\]
Therefore, we find that for all $s\in[t,t+\tau]$
\begin{align*}
  \|S(x_1 (s)) - S(x_2 (s))\| 
  &\le L \| S\big(x_1 - x_2\big)(s)\| \\
  &\le L \gamma(\|\big(\zeta_1 - \zeta_2\big)|_{[-h,s]}\|_\infty)\\
  &\le L_2 \|\big(\zeta_1 - \zeta_2\big)|_{[t,t+\tau]}\|_\infty.
\end{align*}
Furthermore, by linearity and (P1) we have
\begin{align*}
  \|Cx_1(s) - Cx_2(s)\| &= \|(\fh\ast(\zeta_1-\zeta_2))(s)\| \\
  &\le C_{\fh} \|\big(\zeta_1 - \zeta_2\big)|_{[t,t+\tau]}\|_\infty
\end{align*}
for all $s\in[t,t+\tau]$. Now define $\hat c := \tilde c + L_2 + C_{\fh}$ and
\[
    L_3 := \sup \setdef{ \|F'(z)\| }{ \left\| z - \begin{pmatrix} \tilde T (\tilde \xi)(t)\\ S(\tilde x)(t)\\ C\tilde x(t)\end{pmatrix}\right\| \le \hat c \delta}
\]
and set
\[
    c := \hat c L_3.
\]
Then we have
\[
    \|T(\zeta_1)(s) - T(\zeta_2)(s)\| \le c \|\big(\zeta_1 - \zeta_2\big)|_{[t,t+\tau]}\|_\infty
\]
for all $s\in[t,t+\tau]$ and this finishes the proof of the theorem.
\end{proof}

It is a consequence of Theorem~\ref{mainthm} that the operator~$T$ defined in Example~\ref{Ex:TranspEqn} satisfies $T\in\cT_0^{1,1}$. As an additional example, note that it is implicitly shown in~\cite{BergPuch19app} that the operator associated with the internal dynamics of a linearized model of a moving water tank system belongs to the class $\cT_h^{\ell,q}$. In fact, there it is shown that~(P1) is satisfied since the transfer function belongs to the {\it Callier-Desoer} class, cf.\ \cite[Sec.~7.1]{CurtZwar95}.

Concluding this section, we consider a class of \emph{linear} infinite dimensional systems, which can be transformed into a Byrnes-Isidori form, which was introduced in~\cite{IlchSeli16}:
\begin{equation}\label{eq:Abc}
\begin{aligned}
\dot{x}(t)& =Ax(t)+ bu(t),\qquad t\geq 0,  \\
 y(t)& = \langle x(t),c\rangle \,,
\end{aligned}
\end{equation}
where $(A,b,c)$ satisfy, for some $r\in\N$, the assumptions
\begin{enumerate}
\item[{(A1)}] \label{item:A}
$A : \cD(A) \subseteq H \to H $ is  the generator of a $\cC_0$-semigroup~$(T(t))_{t\geq0}$  in a real Hilbert space~$H$ with inner product $\langle\cdot,\cdot\rangle$,
\item[{(A2)}] \label{item:bc}
 $b \in \cD(A^r)$ and $c \in \cD\big((A^*)^r\big)$,
\item[{(A3)}] \label{item:reldeg}
      $\gamma:= \langle A^{r-1}b,c\rangle \ne 0 \  \text{and}  \  \langle A^jb,c\rangle=0$ for all $j=0,1,\ldots,r-2$.
\end{enumerate}
We show that the systems~\eqref{eq:Abc} belong to the class of systems~\eqref{eq:nonlSys}, provided the internal dynamics satisfy a certain BIBO stability assumption. To this end, observe that by~\cite[Thm.~2.6]{IlchSeli16}, system~\eqref{eq:Abc} can be rewritten as
\begin{align*}
  y^{(r)}(t) &= \sum_{i=0}^{r-1} P_i y^{(i)}(t) + S \eta(t) + \gamma u(t),\\
  \dot \eta(t) &= Q \eta(t) + R y(t),\quad \eta(0) = \eta^0,
\end{align*}
where $P_i\in\R$ for $i=0,\ldots,r-1$, $S\in\cL(\hat H;\R)$, $R\in\cL(\R;\hat H)$ and $Q:\cD(Q)\subseteq \hat H \to \hat H$ is the generator of a $\cC_0$-semigroup on~$\hat H$, where $\hat H$  is some real Hilbert space, and $\eta^0\in\cD(Q)$. As a BIBO stability assumption we impose that the transfer function $H(s) = S(sI-Q)^{-1} R$ has inverse Laplace transform which is a Borel measure with bounded total variation.

We may now define the operator~$T$ by
\[
    T(\zeta) := S\eta,\quad \zeta\in\cC(\R_{\ge 0};\R),
\]
where~$\eta$ is the mild solution of $\dot \eta(t) = Q \eta(t) + R \zeta(t)$ with $\eta(0) = \eta^0$. It is clear that~$R$ is a $L^2$-admissible control operator, $S$ is a $L^2$-admissible observation operator and the system~$(Q,R,S)$ is well-posed and regular. Since assumptions~(P2)--(P4) are trivially satisfied in our case, it thus follows that $T\in\cT_0^{1,1}$.

As a consequence, the class of infinite-dimensional systems~\eqref{eq:Abc} is indeed contained in the system class~\eqref{eq:nonlSys}. Moreover, the class of operators $\cT_h^{\ell,q}$ in particular covers operators coming from linear PDE systems as above, but also allows for much more general (and even nonlinear) equations.

%
\section{Simulation}\label{Sec:Sim}
%

We revisit Example~\ref{Ex:TranspEqn} and illustrate our results by a simulation of the funnel controller~\eqref{eq:fun-con} for system~\eqref{eq:sys-sim}. For the simulation we have chosen $\fh(\xi)=f(\xi){\rm d}\xi$ with $f(\xi) = e^{-\xi}/\sqrt{\xi}$, which is integrable but not square integrable on~$\R_{\ge 0}$. Furthermore, we use the parameters $c = \gamma = 1$ and the reference signal
$$y_{\rm ref}(t)=\cos t,\quad t\ge 0.$$
The initial value is chosen as $y(0) = 0$ and for the controller~\eqref{eq:fun-con} we chose the funnel function
\[
  \varphi(t)=\big(2 e^{-2t}+0.1 \big)^{-1},\quad t\ge 0.
\]
Clearly, the initial error lies within the funnel boundaries as required in Theorem~\ref{Thm:FunCon-nonl}. Furthermore, by Theorem~\ref{mainthm} the operator~$T$ satisfies~(N4) and hence funnel control is feasible.

\begin{figure}[h!tb]
\vspace{2mm}
  \centering
  \subfloat[Fig.~\ref{fig:sim}a: Performance funnel with tracking error~$e$ and generated input function~$u$.]
{
\centering
  \includegraphics[width=8.5cm]{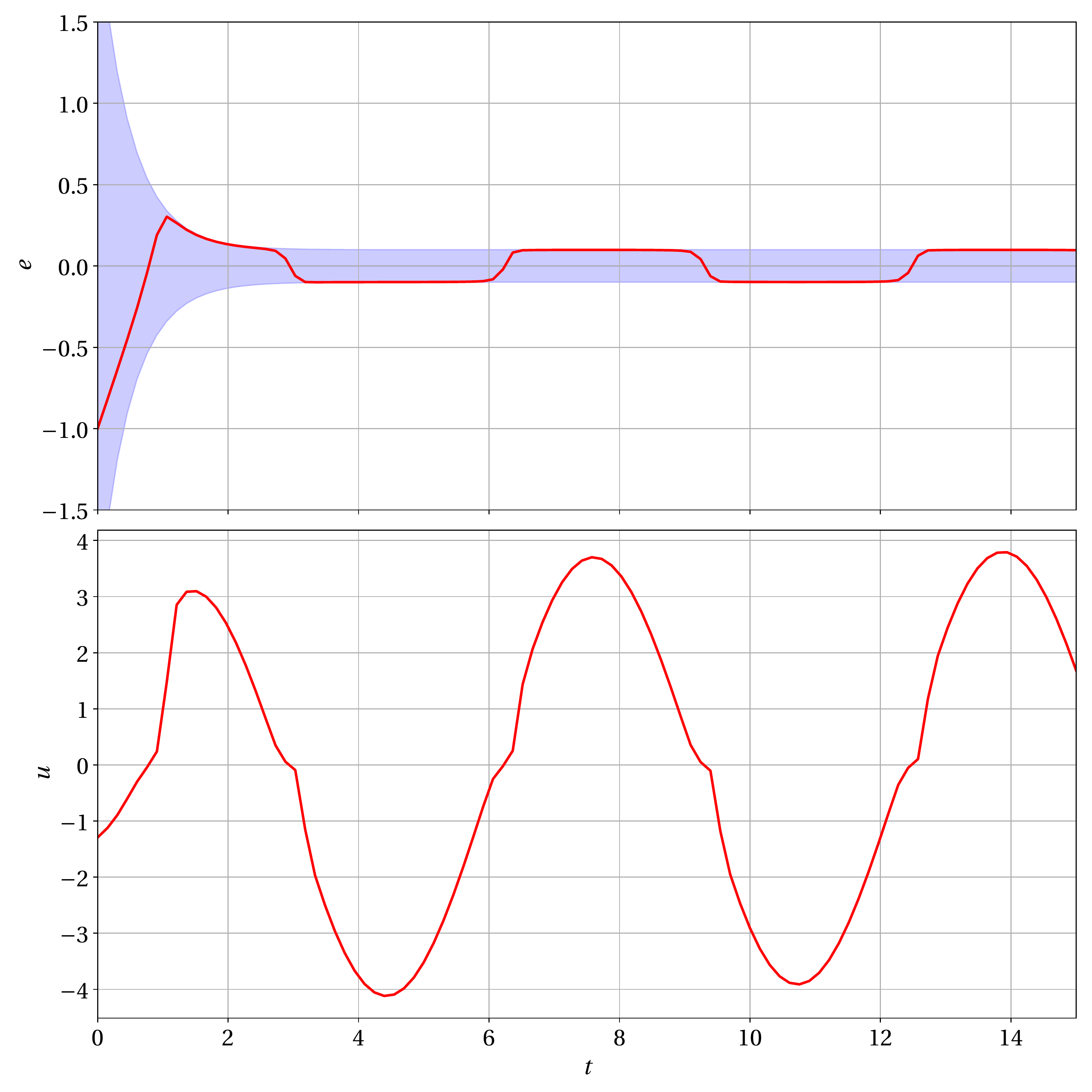}
\label{fig:sim-eu}
}\\
\subfloat[Fig.~\ref{fig:sim}b: State~$z$ of the PDE.]
{
\centering
 \hspace*{-7mm} \includegraphics[width=9.5cm]{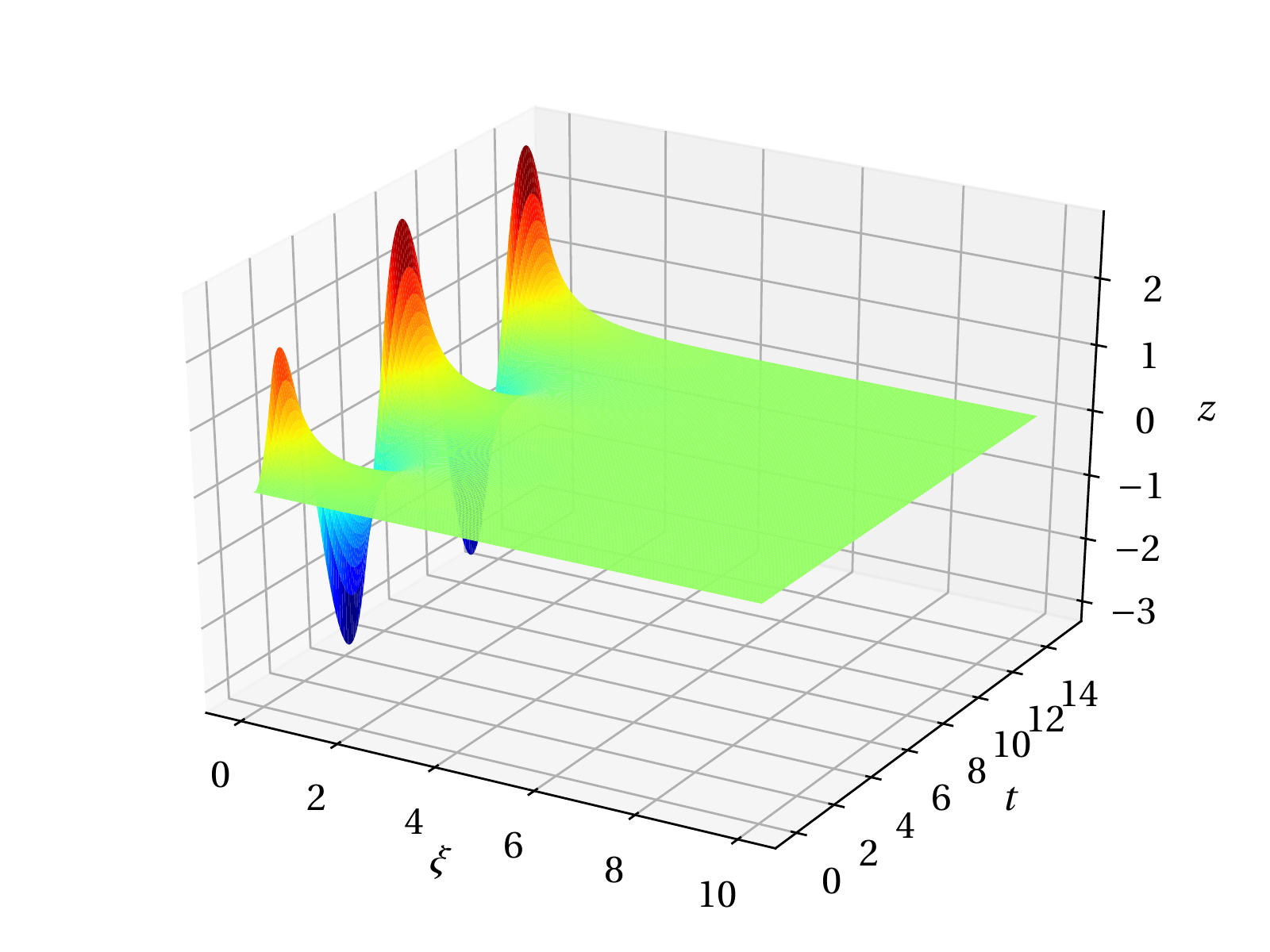}
\label{fig:sim-z}
}
\caption{Simulation of the funnel controller~\eqref{eq:fun-con} for the system~\eqref{eq:sys-sim}.}
\label{fig:sim}
\end{figure}

The PDE is solved using explicit finite differences with a grid in~$t$ with $M=1000$ points for the interval $[0,T]$, where $T=15$, and a grid in~$\xi$ with $N=\lfloor M (b-a)/(\alpha T)\rfloor$ points for $\alpha = 0.4$ and $a=0$, $b=10$. The method has been implemented in Python and the simulation results are shown in Fig.~\ref{fig:sim}.

\captionsetup[subfloat]{labelformat=empty}

It can be seen that even in the presence of infinite-dimensional internal dynamics which are not exponentially stable a prescribed performance of the tracking error can be achieved with the funnel controller~\eqref{eq:fun-con}. At the same time the input generated by the controller is bounded with a very good performance.

%
\section{Conclusion}\label{Sec:Concl}
%

In the present paper we considered the question which classes of systems with infinite-dimensional internal dynamics are encompassed by the abstract system class~\eqref{eq:nonlSys} for which funnel control is feasible by Theorem~\ref{Thm:FunCon-nonl}. We have defined a class of operators~$\cT_h^{\ell,q}$, which model the internal dynamics of the system, that encompass BIBO stable linear and nonlinear PDEs. The corresponding nonlinear observation operators are assumed to satisfy a certain growth bound, while the linear observation operator may be unbounded. For the latter we additionally assumed that the resulting system is regular and well-posed such that the inverse Laplace transform of its transfer function defines a measure with bounded total variation. In Theorem~\ref{mainthm} we have proved that any operator belonging to~$\cT_h^{\ell,q}$ satisfies the conditions of the system class~\eqref{eq:nonlSys}.

Several extensions of the operator class~$\cT_h^{\ell,q}$ and Theorem~\ref{Thm:FunCon-nonl} may be investigated in future research. In particular, extensions to nonlinear PDE systems with unbounded observation operators are of interest as well as systems with infinite-dimensional input and output spaces which do not have an integer-valued relative degree.

\bibliographystyle{elsarticle-harv}

\end{document}